\UseRawInputEncoding
\documentclass[11pt]{amsart}
\usepackage{amssymb,latexsym}
\usepackage{mathrsfs}
\usepackage{amsfonts}
\usepackage{color}
\usepackage{latexsym,amsmath,amssymb}
\usepackage{hyperref}

\newtheorem{theorem}{Theorem}[section]
\newtheorem{lemma}[theorem]{Lemma}
\newtheorem{proposition}[theorem]{Proposition}
\newtheorem{corollary}[theorem]{Corollary}

\theoremstyle{definition}

\theoremstyle{remark}

\numberwithin{equation}{section}
\begin{document}
\title[critical p-Laplace equation]{ On the classification of entire solutions to the critical p-Laplace equation}

\author{Qianzhong Ou}
\address{ School of Mathematics and Statistics\\
          Guangxi Normal University\\
           Guilin, 541006, Guangxi Province, China.}
\email{ouqzh@gxnu.edu.cn}


\maketitle
\begin{abstract}
 Under the assumption of finite energy, positive solutions to the critical p-Laplace equation in $\mathbb{R}^n$ for $1< p<n$ have been classified
 completely by moving plane method. In this paper, the author provide a new approach to obtain the same classification results for $\frac{n+1}{3}\leq p<n$,
 without any further assumptions.

\end{abstract}




\section{Introduction}

\setcounter{equation}{0}

Consider the following  equation

\begin{equation}\label{1.1}
  -\triangle_p u =u^{\alpha-1}          \qquad \text{in} \quad \mathbb{R}^n,
\end{equation}
where  $n\geq 2$, $1<p<n$ and $\triangle_p$  is the usual p-Laplace operator, explicitly

$$\triangle_p u=\verb"div"(|\nabla u|^{p-2}\nabla u),$$
\noindent with $\nabla $ denoted the gradient operator.

The equation  (\ref{1.1})  has been the object of several studies in the differential geometry
and in the PDE's communities. Especially, for $\alpha=\frac{np}{n-p} $,  (\ref{1.1})  is related to the study of the critical points
of the Sobolev inequality (see e.g. the survey \cite{Ro2022}) and  $p^{*}:=\frac{np}{n-p}$ is called the Sobolev exponent,
meanwhile for $p = 2$, to the Yamabe problem (see e.g. the survey \cite{LP1987}).
Indeed, the p-Laplacian operator also appears in the studies of physics (see e.g. \cite{DF1994} and the references therein)
and  stochastic models (see e.g. \cite{SZi2022}).

In this paper, we focus on the classifications of positive solutions to  (\ref{1.1}),
especially for the critical case $\alpha=p^{*}$. It is well known that such issue is crucial in many applications
such as {\it a priori} estimates, blow-up analysis and asymptotic analysis. Let us recall some known results on this issue in below.

For the subcritical case $1<\alpha<p^{*}$, the nonnegative solutions to (\ref{1.1})
have been characterized completely by two splendid papers: Gidas-Spruck \cite{GS1981} for $p=2$ and Serrin-Zou \cite{SZ2002} for general $1<p<n$,
where they showed that all the nonnegative (weak) solutions to (\ref{1.1}) are trivial.
While for the critical case $\alpha=p^{*}$, there are nontrivial 2-parameters family of solutions of  (\ref{1.1}) as follows

\begin{equation}\label{1.2}
 U_{\lambda,x_0}(x):=\Big( \frac{\lambda^{\frac{1}{p-1}}n^{\frac{1}{p}}\big(\frac{n-p}{p-1}\big)^{\frac{p-1}{p}}}{\lambda^{\frac{p}{p-1}}+|x-x_0|^{\frac{p}{p-1}}}\Big)^{\frac{n-p}{p}}     \qquad \text{with} \,\,\lambda>0,\, x_0\in \mathbb{R}^n.
\end{equation}

Now we focus on the classifications of positive solutions to the critical case $\alpha=p^{*}$, i.e., the following critical p-laplace equation
\begin{equation}\label{1-1}
  -\triangle_p u =u^{p^*-1}          \qquad \text{in} \quad \mathbb{R}^n
\end{equation}
for $1<p<n$. Then an interesting and challenging problem is:

{\it whether any  positive solution of  (\ref{1-1}) must be of the form  (\ref{1.2})?}

In case of $p=2$, this problem was solved by Caffarelli-Gidas-Spruck \cite{CGS89} via the method of moving planes
and the Kelvin transform (see also \cite{CL1991, GNN1979, LZ2003, LZhu1995, Ob1971}).

For $p\neq 2$, the problem is quasilinear and the Kelvin transform is not available, which makes it more complicated then the semilinear case.
Then, under the additional assumption of finite energy, the moving plane method were also exploited  by B. Sciunzi \cite{Sc2016} and J. V$\acute{e}$tois \cite{Ve2016}
 to show that any positive  weak solution of (\ref{1-1})  must be of the form $u(x)=U_{\lambda,x_0}(x)$.
 Recall that  the energy associated to  (\ref{1-1})  is given by

 $$E_{\mathbb{R}^n}(u):= \frac{1}{p}\int_{\mathbb{R}^n}|\nabla u|^p+ \frac{1}{p^{*}}\int_{\mathbb{R}^n} u^{p^{*}},$$
\noindent and the equation (\ref{1-1}) is just the Euler-Lagrange equation associated to this energy functional.

The same result of \cite{Sc2016, Ve2016}  was extended by  Ciraolo-Figalli-Roncoroni \cite{CFR2020}  to the anisotropic setting and in convex cones of $\mathbb{R}^n$.
Note that the method used in \cite{CFR2020} is different than that in \cite{Sc2016, Ve2016} and is closed to that in \cite{GS1981, SZ2002}.
In fact, the method used in  \cite{GS1981, SZ2002} may be originally due to Obata \cite{Ob1971} when he studied the rigidity of the conformal metrics on the sphere.
Roughly speaking, the soul of this method is to construct some suitable vector fields and
then deduce an integral estimate for these vector fields via integration by parts with appropriate test functions.
The same method had been also used successfully in the analogous problems on Heisenberg group (see e.g. \cite{JL1988, MO2020})
and in some fully nonlinear problems  (see e.g. \cite{Gon06, OU2013}).

But for $p\neq 2$ and without any further assumptions, the problem is still completely open, until recently Catino-Monticelli-Roncoroni \cite{CMR2022}
gave a result in dimensions $n=2,3$ for $\frac{n}{2}<p<2$. In \cite{CMR2022}, they also deal with this problem with the solutions satisfying
suitable conditions at infinity. The argument in  \cite{CMR2022} was borrowed from Serrin-Zou \cite{SZ2002} and Ciraolo-Figalli-Roncoroni \cite{CFR2020},
especially they used nearly the same fundamental integral inequality as in \cite{SZ2002} (see proposition 2.2 in \cite{CMR2022}).

   In this paper, we deal with the problem for all dimensions  $n\geq 2$ and a wider range of $p$, say $\frac{n+1}{3}\leq p<n$.
Our strategy is adopt the main idea of \cite{Ob1971,GS1981,SZ2002,CFR2020} and \cite{CMR2022},
 but we will adopt a new test function and generalize the fundamental integral inequality appeared in \cite{SZ2002, CMR2022},
 and then deal with the ``error" terms by more careful analysis. Our main result is as follows.

\begin{theorem}\label{Thm1}
Let $1<p<n$ and $u$ be a positive weak solution of (\ref{1-1}). Assume in addition $\frac{n+1}{3}\leq p<n$,
 then $u$ must be of the form  (\ref{1.2}), i.e., $u=U_{\lambda,x_0}$ for some $\lambda>0$ and $x_0\in \mathbb{R}^n$.
\end{theorem}

Note that a function $u\in W^{1,p}_{loc}(\mathbb{R}^n)\cap L^{\infty}_{loc}(\mathbb{R}^n)$ is said to be a weak solution of  (\ref{1-1})  if

\begin{equation}\label{1.3}
  \int_{\mathbb{R}^n}|\nabla u|^{p-2}\langle \nabla u,\nabla \psi\rangle-  \int_{\mathbb{R}^n} u^{p^{*}-1}\psi =0 \qquad \forall\, \psi \in W^{1,p}_{0}(\mathbb{R}^n).
\end{equation}

Here we mention some well-known facts about solutions of (\ref{1-1}), that are instrumental in the proof of theorem \ref{Thm1}.
By the strong maximum principle, all nonnegative nontrivial solutions of (\ref{1-1}) must be strictly positive.
So in what follows we shall always be concerned with positive weak solutions. For a positive weak solution $u$ of (\ref{1-1}), we have

\begin{equation}\label{1.4}
  u(x) \geq C |x|^{-\frac{n-p}{p-1}} \qquad \text{for}  \quad |x|>1,
\end{equation}
where the positive constant $C$ depending only on $n,p $ and $\min_{|x|=1}u$.
In fact, the estimate  (\ref{1.4}) had been derived for positive weak sub-p-harmonic functions (see lemma 2.3 in \cite{SZ2002}).
For the regularity of a positive weak solution  $u$ of (\ref{1-1}) (see e.g. \cite{ACF2022}), we have

\begin{equation}\label{1.5}
  u\in W^{2,2}_{loc}(\mathbb{R}^n)\cap C^{1,\theta}_{loc}(\mathbb{R}^n),
\end{equation}
for some $\theta\in(0,1)$ and, in addition,

\begin{equation}\label{1.6}
  |\nabla u|^{p-2}\nabla u\in W^{1,2}_{loc}(\mathbb{R}^n),
\end{equation}
and

\begin{equation}\label{1.7}
 |\nabla u|^{p-2}\nabla^2u\in L^{2}_{loc}(\mathbb{R}^n).
\end{equation}
Moreover, if we set

$$\Omega_{cr}=\{x\in \mathbb{R}^n|\, \nabla u(x)=0\},$$
\noindent then $\Omega_{cr}$ has zero measure and any weak solution is actually $C^{\infty}$  on $\Omega^{c}_{cr}$ by the bootstrap argument.

The paper is organized as follows. In section 2, we will give some preparation calculations. We shall introduce
some vector fields and present some divergent identities and some inequalities on these vector fields,
 especially, we give the fundamental integral inequality.
Then, using this key integral inequality combining with previous inequalities on the vector fields,
we will prove theorem \ref{Thm1} in section 3 by careful analysis on the ``error" terms.

\section{Preliminaries}

\setcounter{equation}{0}
\setcounter{theorem}{0}

There are two parts in this section.
In the first part, we will take a transformation for any positive weak solution of  (\ref{1-1}).
Then we introduce some vector fields associated to the new transformed function and give a fundamental integral estimate for these vector fields.
In the second part, we proof some inequalities on those vector fields, which will be used to estimate the ``error" terms in the key integral estimate.

\subsection{Vector fields and the fundamental integral identity}
 Let $u>0$ be any weak solution of (\ref{1-1}) and take $v=u^{-\frac{p}{n-p}}$ . Then $v$ satisfies, also in the weak sense,

\begin{equation}\label{2.1}
  \triangle_p v =\, \frac{p-1}{p}n v^{-1}|\nabla v|^p +(\frac{p}{n-p})^{p-1}v^{-1}        \qquad\quad \text{in} \quad \mathbb{R}^n.
\end{equation}

Clearly $v$ inherits some properties from $u$. In particular $v$ share the same critical set $\Omega_{cr}$ with $u$ and $v\in C^{\infty}(\Omega^c_{cr})$.
Moreover, by (\ref{1.4}) we have

\begin{equation}\label{2.2}
  v(x)\leq C|x|^{\frac{p}{p-1}} \qquad \text{for} \quad |x|>1,
\end{equation}
with the constant $C$ depending only on $n,p$ and  $\max_{|x|=1}v(x)$ and by (\ref{1.5})-(\ref{1.7}) we also have

\begin{equation}\label{2.3}
  v\in W^{2,2}_{loc}(\mathbb{R}^n)\cap C^{1,\theta}_{loc}(\mathbb{R}^n),
\end{equation}
for some $\theta\in(0,1)$ and

\begin{equation}\label{2.4}
  |\nabla v|^{p-2}\nabla v\in W^{1,2}_{loc}(\mathbb{R}^n),
\end{equation}

\begin{equation}\label{2.5}
 |\nabla v|^{p-2}\nabla^2v\in L^{2}_{loc}(\mathbb{R}^n).
\end{equation}

Now we introduce the following vector fields

$$X^i=|\nabla v|^{p-2} v_i,$$

$$E_{ij}=X^{i}_{,j}-\frac{1}{n}X^{k}_{,k}\delta_{ij},\qquad E_j=v^{-1}v_iE_{ij},$$
\noindent where and in the sequel, we adopt the Einstein convention of summation over repeated indices.
Since $p>1$,  it is understood in the usual way that $X^i$ and then $E_{ij}$ is identically zero in $\Omega_{cr}$.
Obviously, the matrix $E=\{E_{ij}\}$ is trace free, i.e., $\textbf{Tr}E=E_{ii}\equiv 0$, but may not be symmetric
 (this differs from the case $p=2$ and hence we must be more careful to deal with the ``error" term).

Denote the function

 $$g=av^{-1}|\nabla v|^p +bv^{-1}$$

\noindent with $a=\frac{p-1}{p}n ,\, b=(\frac{p}{n-p})^{p-1}$. Then we have the equation

$$ \triangle_p v= X^{k}_{,k} =g             \qquad \text{in} \quad \mathbb{R}^n$$
\noindent in the weak sense, that is

\begin{equation}\label{2.6}
   -\int_{\mathbb{R}^n} X^{k}\psi_k=  \int_{\mathbb{R}^n} g\psi  \qquad \forall\, \psi \in W^{1,p}_{0}(\mathbb{R}^n).
\end{equation}
Also we have
$$E_{ij}=X^{i}_{,j}-\frac{1}{n}g\delta_{ij},\qquad E_j=v^{-1}v_iX^{i}_{,j}-\frac{1}{n}gv^{-1}v_j.$$

 Since $v\in C^{\infty}(\Omega^c_{cr})$, so before presenting our key integral estimate,
 we give some differential identities on $\Omega^c_{cr}$ and these may be helpful to understand the key integral estimate itself.

\begin{lemma}\label{lem-1}
With the notations as in above, then in $\Omega^c_{cr}$ we have

\,\,(i)\,\,\,$g_j=nv^{-1}v_iE_{ij}=nE_j$;

\,(ii)\,\,$E_{ij,i}=\frac{n-1}{n}g_j=(n-1)E_j$;

(iii)\,$(X^jE_{ij})_{,i}=E_{ij}E_{ji}+(n-1)X^jE_j=\mathbf{Tr}\{E^2\}+(n-1)X^jE_j$.
\end{lemma}

\vspace{20pt}
$\mathbf{Proof\,\, of\,\, Lemma\,\, \ref{lem-1}}$ \qquad First we have

\begin{equation}\label{2.7}
\begin{split}
  v^{-1}v_iE_{ij}=& v^{-1}v_i\big(X^{i}_{,j}-\frac{1}{n}g\delta_{ij}\big)\\
                 =& v^{-1}v_i\big((p-2)|\nabla v|^{p-4}v_kv_{kj}v_i +|\nabla v|^{p-2}v_{ij}-\frac{1}{n}g\delta_{ij}\big)\\
                 =& (p-1) v^{-1}|\nabla v|^{p-2}v_iv_{ij}-\frac{1}{n}gv^{-1}v_{j}.
\end{split}
\end{equation}
On the other hand,

\begin{equation}\label{2.8}
\begin{split}
   g_j=& \big( a|\nabla v|^pv^{-1} +bv^{-1} \big)_{j}\\
      =& \big( a|\nabla v|^p +b \big)_{j}v^{-1}+\big( a|\nabla v|^p +b \big)(v^{-1})_{j}\\
      =& ap|\nabla v|^{p-2}v_iv_{ij}v^{-1}+\big( a|\nabla v|^p +b \big)(-v^{-2}v_{j})\\
      =& n(p-1)|\nabla v|^{p-2}v_iv_{ij}v^{-1}- gv^{-1}v_{j}.
\end{split}
\end{equation}
Comparing (\ref{2.7}) with (\ref{2.8}) we get (i).

Using the equation (\ref{2.2}) and (i)  we can prove (ii) as follows

\begin{equation}\label{2.9}
\begin{split}
  E_{ij,i}=& (X^i_{,j}-\frac{1}{n}g\delta_{ij})_{,i}\\
          =& X^i_{,ji}-\frac{1}{n}g_j\\
          =& X^i_{,ij}-\frac{1}{n}g_j\\
          =& g_j-\frac{1}{n}g_j\\
          =& (n-1)E_j.
\end{split}
\end{equation}

Using (i),(ii) and notice that $E$ is trace free, we get further

\begin{equation}\label{2.10}
\begin{split}
(X^jE_{ij})_{,i}=& X^j_{,i}E_{ij} +X^jE_{ij,i}\\
          =& E_{ji}E_{ij}+(n-1)X^jE_j.
\end{split}
\end{equation}
Thus (iii) is also valid.\qed

Furthermore, also on $\Omega^c_{cr}$ we have

\begin{lemma}\label{lem-3}
\begin{equation}\label{2.11}
\big( v^{q}g^{m}X^j\big)_{,j}=(a+q)v^{q-1}g^m|\nabla v|^p+bv^{q-1}g^{m}+nmv^{q}g^{m-1}X^{i}E_i
\end{equation}
and

\begin{equation}\label{2.12}
\big( v^{q}g^{m}X^jE_{ij}\big)_{,i}= v^{q}g^m \mathbf{Tr}\{E^2\} +nmv^{q}g^{m-1}X^{j}E_{ij}E_i +(n-1+q)v^{q}g^{m}X^{j}E_j,
\end{equation}
where $q$, $m$ are constants underdetermined.
\end{lemma}

\vspace{20pt}
$\mathbf{Proof\,\, of\,\, Lemma\,\, \ref{lem-3}}$ \qquad Using (i) in lemma \ref{lem-1}
and the equation (\ref{2.2}) we deduce  (\ref{2.11}) as follows

\begin{equation}\label{2.13}
\begin{split}
  \big(v^{q}g^mX^{i}\big)_{,i}
  =& qv^{q-1}g^mX^{i}v_i+mv^{q}g^{m-1}g_iX^{i}+v^{q}g^mX^{i}_{,i}\\
  =& qv^{q-1}g^m|\nabla v|^p+mv^{q}g^{m-1}(nE_i)X^{i}+v^{q}g^{m+1}\\
  =& (a+q)v^{q-1}g^m|\nabla v|^p+bv^{q-1}g^{m}+nmv^{q}g^{m-1}X^{i}E_i.
\end{split}
\end{equation}
By (i),(iii) in lemma \ref{lem-1} we also obtain

\begin{equation}\label{2.14}
\begin{split}
 \,& \big(v^{q}g^mX^{j}E_{ij}\big)_{,i}\\
  =& v^{q}g^m\big(X^{j}E_{ij}\big)_{,i} + \big(v^{q}g^m\big)_{,i}X^{j}E_{ij}\\
  =& v^{q}g^m\big( \mathbf{Tr}\{E^2\}+(n-1)X^jE_j\big) + \big(qv^{q-1}v_ig^m+mv^qg^{m-1}g_i\big) X^{j}E_{ij}\\
  =& v^{q}g^m\big( \mathbf{Tr}\{E^2\}+(n-1)X^jE_j\big) + qv^{q}g^mX^{j}E_{j} +mv^qg^{m-1}(nE_i ) X^{j}E_{ij}\\
  =& v^{q}g^m\mathbf{Tr}\{E^2\} +nmv^{q}g^{m-1}X^{j}E_{ij}E_i +(n-1+q)v^{q}g^{m}X^{j}E_j.
\end{split}
\end{equation}
This is  (\ref{2.12}). \qed

Next, along these lines of \cite{Ob1971,GS1981,SZ2002,CFR2020,CMR2022}, the idea is to apply the lemma \ref{lem-3}
(especially  (\ref{2.12}) with $q=1-n$) and integrate the identities over $\mathbb{R}^n$ after multiplying suitable
test functions. Due to the lack of regularity of $v$, lemma \ref{lem-3} cannot be
applied directly but we can still prove its integral counterpart. First we have the following fundamental integral inequality
which is a generalization of that in \cite{CMR2022, SZ2002} (see proposition 2.2 in \cite{CMR2022} or proposition 6.2 in \cite{SZ2002} ).

\begin{proposition}\label{pro2.4}
Let $u$ be any positive weak solution of (\ref{1-1}) and using the notations as before, then for every $0\leq \varphi\in C^{\infty}_0(\mathbb{R}^n)$ we have
\begin{equation}\label{2.15}
\int_{\mathbb{R}^n}\varphi v^{1-n}g^m \mathbf{Tr}\{E^2\}+nm\int_{\mathbb{R}^n}\varphi v^{1-n}g^{m-1}X^{j}E_{ij}E_i\leq -\int_{\mathbb{R}^n} v^{1-n}g^{m}X^jE_{ij}\varphi_i .
\end{equation}
\end{proposition}

\vspace{20pt}
$\mathbf{Proof\,\, of\,\, Proposition\,\, \ref{pro2.4}}$ \qquad In case $m=0$ the result follows from
proposition 2.2 in \cite{CMR2022} (see also proposition 6.2 in \cite{SZ2002}), that is

\begin{equation}\label{2.16}
 \int_{\mathbb{R}^n}\varphi v^{1-n}\mathbf{Tr}\{E^2\} \leq -\int_{\mathbb{R}^n} v^{1-n}X^jE_{ij}\varphi_i .
\end{equation}

For $m\neq 0$ we argue by approximation.
So, for $\epsilon>0$ we define $v^{\epsilon}=v\ast \rho^{\epsilon}$, where $\rho^{\epsilon}$ is a standard mollifier.
We also denote

 $$g^{\epsilon}=a(v^{\epsilon})^{-1}|\nabla v^{\epsilon}|^p +b(v^{\epsilon})^{-1}.$$

Now replacing $\varphi$ with $(g^{\epsilon})^m\varphi$ in (\ref{2.16}) we have

\begin{equation}\label{2.17}
 \int_{\mathbb{R}^n}(g^{\epsilon})^m\varphi v^{1-n}\mathbf{Tr}\{E^2\}
 \leq -\int_{\mathbb{R}^n} v^{1-n}X^jE_{ij}\big[m(g^{\epsilon})^{m-1}(g^{\epsilon})_i\varphi+(g^{\epsilon})^m\varphi_i\big].
\end{equation}

Notice that $E_{ij}\in L^{2}_{loc}(\mathbb{R}^n)$ and $g\in W^{1,2}_{loc}(\mathbb{R}^n)$,
we have $g^{\epsilon}\rightarrow g$ and $(g^{\epsilon})_i\rightarrow g_i=nv^{-1}v_jE_{ji}=nE_i$ in $L^{2}_{loc}(\mathbb{R}^n)$ as  $\epsilon\rightarrow 0$.
 Then we obtain (\ref{2.15}) from (\ref{2.17}) easily by letting $\epsilon\rightarrow 0$.\qed

The following is the integral counterpart of (\ref{2.11}) with $m=0$,
which will be needed in next section to deal with the ``error" term on the right hand side of (\ref{2.15}).

\begin{lemma}\label{lem5}
\begin{equation}\label{2.18}
(a+1-q)\int_{\mathbb{R}^n}v^{-q}|\nabla v|^p\psi +b\int_{\mathbb{R}^n} v^{-q}\psi = -\int_{\mathbb{R}^n} v^{1-q} X^j\psi_j \qquad \forall\, \psi \in W^{1,p}_{0}(\mathbb{R}^n).
\end{equation}
\end{lemma}

\vspace{20pt}
$\mathbf{Proof\,\, of\,\, Lemma\,\, \ref{lem5}}$ \qquad Replacing $\psi$ with $v^{1-q}\psi$
in  (\ref{2.6}) and then an elementary computation concludes the result.\qed

\subsection{Inequalities on Vector fields}

In this subsection, we will prove some inequalities on those vector fields introduced  in the last subsection,
which will be used to estimate the ``error" term in the key integral inequality  (\ref{2.15}).

First we give an inequality for some general matrices.

\begin{lemma}\label{lem-6}
Let $A=\{a_{ij}\}$, $B=\{b_{ij}\}$, $C=\{c_{ij}\}$ be $n\times n$ square matrices. If $A$ is positive definite and  diagonal, precisely,
if  $a_{ij}=\lambda_i\delta_{ij}$ with $\lambda_i>0,\, i=1,2,\cdots,n$, then we have

\begin{equation}\label{2.19}
2\mathbf{Tr}\{BAC\}\leq \frac{\Lambda^2}{\lambda^2}\mathbf{Tr}\{BB^{t}\}+\lambda^2\mathbf{Tr}\{CC^t\},
\end{equation}
where $\Lambda=\max\{\lambda_1,\lambda_2,\cdots,\lambda_n\}$, $\lambda=\min\{\lambda_1,\lambda_2,\cdots,\lambda_n\}$.
\end{lemma}

\vspace{20pt}
$\mathbf{Proof\,\, of\,\, Lemma\,\, \ref{lem-6}}$ \qquad First we have

\begin{equation}\label{2.20}
\mathbf{Tr}\{BAC\} = \sum^n_{i,j,k=1} b_{ij}a_{jk}c_{ki}
 = \sum^n_{i,j,k=1}b_{ij}\lambda_j\delta_{jk}c_{ki}= \sum^n_{i,j =1}b_{ij}\lambda_jc_{ji}.
\end{equation}
The Cauchy-Schwarz inequality shows

$$
 2\sum^n_{i,j =1}b_{ij}\lambda_jc_{ji}
\leq  \sum^n_{i,j=1} \big(\frac{\Lambda}{\lambda}b_{ij}\big)^2
                        +\sum^n_{i,j=1}\big(\frac{\lambda}{\Lambda}\lambda_j c_{ji}\big)^2.
$$
\noindent Therefor

\begin{equation}\label{2.21}
\begin{split}
2\mathbf{Tr}\{BAC\}
\leq&\, \sum^n_{i,j=1} \big(\frac{\Lambda}{\lambda}b_{ij}\big)^2
                        +\sum^n_{i,j=1}\big(\frac{\lambda}{\Lambda}\lambda_j c_{ji}\big)^2\\
 =&\, \frac{\Lambda^2}{\lambda^2}\sum^n_{i,j=1} b_{ij}^2+\lambda^2\sum^n_{i,j=1} \frac{\lambda^2_j}{\Lambda^2} c^2_{ji}\\
\leq&\, \frac{\Lambda^2}{\lambda^2}\sum^n_{i,j=1} b_{ij}^2+\lambda^2\sum^n_{i,j=1} c^2_{ji}\\
 =& \frac{\Lambda^2}{\lambda^2}\mathbf{Tr}\{BB^{t}\}+\lambda^2\mathbf{Tr}\{CC^t\} . \qed
\end{split}
\end{equation}

\vspace{20pt}
As an application of lemma \ref{lem-6}, we have the following

\begin{corollary}\label{cor-1}
Let $B$ be any $n\times n$ square matrix, then we have

\begin{equation}\label{2.22}
\mathbf{Tr}\{BE\}\leq c(p)\mathbf{Tr}\{BB^{t}\}+ \mathbf{Tr}\{E^2\},
\end{equation}
where  $c(p)$ is a positive constant depending only on $p$.
\end{corollary}

\vspace{20pt}
$\mathbf{Proof\,\, of\,\, Corollary\,\, \ref{cor-1}}$ \qquad On $\Omega_{cr}$, since $E=0$,  the conclusion is obvious.
In $\Omega^{c}_{cr}$, first we observe $X^i_{,j}=|\nabla v|^{p-2}AH$ with $H=\nabla^2v$ (
the Hessian of $v$ ) and $A=I+(p-2)\frac{\nabla v\bigotimes \nabla v}{|\nabla v|^2}$
being positive definite with eigenvalues  $\lambda_1=p-1, \lambda_i=1, i=2,\cdots,n$.
Then we can rewrite
$E=|\nabla v|^{p-2}AH-\frac{1}{n}gI=AC$ with $C= |\nabla v|^{p-2}H-\frac{1}{n}gA^{-1}$.

By rotating the coordinate system we may assume $\nabla v=(v_1,0,\cdots,0)$. Then $A$ is diagonal, precisely,
the entries $a_{ij}=\lambda_i\delta_{ij}$. Clearly $C$ is symmetric and then we have

 \begin{equation}\label{2.23}
 \begin{split}
\mathbf{Tr}\{E^2\}=&\mathbf{Tr}\{ACAC\}\\
=& \sum^n_{i,j,k,l=1}a_{ij}c_{jk}a_{kl}c_{li}\\
=&\sum^n_{i,j,k,l=1}\lambda_i\delta_{ij}c_{jk}\lambda_k\delta_{kl}c_{li}\\
=&\sum^n_{i,k=1}\lambda_i\lambda_kc^2_{ik}\\
\geq & \lambda^2\mathbf{Tr}\{CC^t\} ,
\end{split}
\end{equation}
where $\lambda=\min\{p-1,1\}$.

On the other hand, with $\Lambda=\max\{p-1,1\}$, using  (\ref{2.19}) we get

\begin{equation}\label{2.24}
\begin{split}
\mathbf{Tr}\{BE\}
 =&\mathbf{Tr}\{BAC\}\\
 \leq& \frac{\Lambda^2}{2\lambda^2}\mathbf{Tr}\{BB^{t}\}+\frac{\lambda^2}{2}\mathbf{Tr}\{CC^t\}.
\end{split}
\end{equation}
 Combining  (\ref{2.24}) with (\ref{2.23}) we obtain (\ref{2.22}).\qed

For the matrix $E$, we also have the following inequality.

\begin{lemma}\label{lem-8}
 \begin{equation}\label{2.25}
\mathbf{Tr}\{E^2\}=\sum^n_{i,j=1}E_{ij}E_{ji}\geq \sum^n_{i,j,k=1}\frac{v_j}{|\nabla v|}E_{ij}\frac{v_k}{|\nabla v|}E_{ki}.
\end{equation}
\end{lemma}

\vspace{20pt}
$\mathbf{Proof\,\, of\,\, Lemma\,\, \ref{lem-8}}$ \qquad Using the computations in  (\ref{2.23}) one has

\begin{equation}\label{2.26}
\mathbf{Tr}\{E^2\}=\sum^n_{i,k=1}\lambda_i\lambda_kc^2_{ik},
\end{equation}
and similarly

\begin{equation}\label{2.27}
\sum^n_{i,j,k=1}\frac{v_j}{|\nabla v|}E_{ij}\frac{v_k}{|\nabla v|}E_{ki}=\sum^n_{k=1}\lambda_1\lambda_kc^2_{1k}.
\end{equation}
Therefor

\begin{equation}\label{2.28}
\mathbf{Tr}\{E^2\}-\sum^n_{i,j,k=1}\frac{v_j}{|\nabla v|}E_{ij}\frac{v_k}{|\nabla v|}E_{ki}=\sum^n_{i=2}\sum^n_{k=1}\lambda_i\lambda_kc^2_{ik}\geq 0. \qed
\end{equation}

\section{Proof of Theorem \ref{Thm1} }

\setcounter{equation}{0}
\setcounter{theorem}{0}

In this section, we will prove theorem \ref{Thm1} by deducing the integral estimate (\ref{3.32}).
 Let $u>0$ be any weak solution of (\ref{1-1}) and take $v=u^{-\frac{p}{n-p}}$.
We will use the notations and the results presented in section 2. In fact, as in \cite{CFR2020, CMR2022}, to show $u=U_{\lambda,x_0}$,
we need only to show $E= 0$. To do this, we will deduce (\ref{3.33}), and then we must have  $E= 0$, since $\mathbf{Tr}\{E^2\}\geq 0$
and the ``=" happens if and only if $E= 0$ (see (\ref{2.23}) or \cite{CFR2020, SZ2002} for more details).

 Let $\eta$ be smooth cut-off functions satisfying:
\begin{equation}\label{3.1}
\begin{cases}
                          \eta\equiv 1  &\verb"in" \,\,B_R,\\
                        0\leq\eta\leq1  &\verb"in" \,\,B_{2R},\\
                          \eta\equiv 0  &\verb"in" \,\,\mathbb{R}^n\backslash B_{2R},\\
     |\nabla \eta|\lesssim \frac{1}{R}  &\verb"in" \,\,\mathbb{R}^n,
\end{cases}
\end{equation}

\noindent where and in the sequel, $B_R$ denotes a
ball in $\mathbb{R}^n$ centered at the origin with radius $R$ ;
and we use ``$\lesssim $" , ``$\backsimeq$"  to replace ``$\leq $", ``$=$",  etc., to drop out some
positive constants independent of $R$ and $v$.

First, we deduce the following integral estimates(see also lemma 2.4 in \cite{SZ2002} for similar results),
which is needed in the proof of theorem \ref{Thm1}.

\begin{lemma}\label{lem3.1}
For $p \leq q<a+1$,

\begin{equation}\label{3.2}
\int_{B_R}v^{-q }|\nabla v|^p\lesssim R^{n-q }
\end{equation}
and for $0<q\leq a+1$

\begin{equation}\label{3.3}
\int_{B_R}v^{-q } \lesssim R^{n-q }.
\end{equation}
\end{lemma}

$\mathbf{Proof\,\, of\,\, Lemma\,\, \ref{lem3.1}}$ \qquad Let $\theta>0$ be a constant big enough and
 $\eta$ be cut off functions as in (\ref{3.1}). Using (\ref{2.18}) with $\psi=\eta^\theta$   we have

\begin{equation}\label{3.5}
\begin{split}
(a+1-q)\int \eta^\theta v^{-q } |\nabla v|^p+b\int \eta^\theta v^{-q }
       = & -\theta\int \eta^{\theta-1}\eta_j  v^{1-q} X^j \\
\lesssim & \frac{1}{R}\int \eta^{\theta-1}v^{1-q}|\nabla v|^{p-1}
\end{split}
\end{equation}
since $|\eta_j X^j|\leq |\nabla\eta| |\nabla v|^{p-1}\lesssim \frac{1}{R}|\nabla v|^{p-1}$.
The  Young's inequality with exponent pair $(\frac{p}{p-1},p)$ shows

\begin{equation}\label{3.6}
\begin{split}
\frac{1}{R}\eta^{\theta-1}v^{1-q}|\nabla v|^{p-1}
   = & \eta^{\theta}v^{-q }\big(\varepsilon^{\frac{p-1}{p}}|\nabla v|^{p-1}\cdot \frac{1}{\varepsilon^{\frac{p-1}{p}}R}\eta^{ -1}v\big) \\
\leq & \eta^{\theta}v^{-q }\big(\varepsilon |\nabla v|^{p}+ \frac{1}{\varepsilon^{p-1} R^p}\eta^{ -p}v^p\big) \\
   = & \varepsilon\eta^{\theta}v^{-q } |\nabla v|^{p}+ \frac{1}{\varepsilon^{p-1} R^{p}}\eta^{ \theta-p}v^{p-q },
\end{split}
\end{equation}
where  $\varepsilon>0$ is a small constant. Therefor

\begin{equation}\label{3.7}
\begin{split}
\,& (a+1-q)\int \eta^\theta v^{-q } |\nabla v|^p+b\int \eta^\theta v^{-q }\\
\lesssim & \varepsilon\int \eta^{\theta}v^{-q } |\nabla v|^{p}+ \frac{1}{\varepsilon^{p-1} R^p}\int \eta^{ \theta-p}v^{p-q }.
\end{split}
\end{equation}
Since $p <q$ implies $\frac{q }{q-p}>1$, then similarly, for the last term in  (\ref{3.7}),
using the Young's inequality with exponent pair $(\frac{q }{q-p}, \frac{q }{p})$ we have

\begin{equation}\label{3.8}
\frac{1}{\varepsilon^{p-1} R^p}\int \eta^{ \theta-p}v^{p-q }
 \leq  \varepsilon\int \eta^\theta v^{-q }+\frac{1}{\varepsilon^{q-1} R^{q }}\int \eta^{ \theta-q }.
\end{equation}
Inserting this into  (\ref{3.7}) yields

\begin{equation}\label{3.9}
\begin{split}
\,& (a+1-q)\int \eta^\theta v^{-q } |\nabla v|^p+b\int \eta^\theta v^{-q }\\
\lesssim & \varepsilon\int \eta^{\theta}v^{-q } |\nabla v|^{p}
          + \varepsilon\int \eta^\theta v^{-q }+\frac{1}{\varepsilon^{q-1} R^{q }}\int \eta^{ \theta-q }.
\end{split}
\end{equation}
Note that (\ref{3.9}) also valid for $q=p $, since then it comes from (\ref{3.7}) immediately.
Recall the definition of $\eta$ in (\ref{3.1}) and if $a+1-q>0 $,
taking $\varepsilon>0$ small enough and $\theta>q $, we obtain from   (\ref{3.9})

\begin{equation}\label{3.10}
\int_{B_R} v^{-q } |\nabla v|^p+ \int_{B_R} v^{-q }
\lesssim  R^{n-q }.
\end{equation}
This implies (\ref{3.2}) and (\ref{3.3}) for $p \leq q<a+1$. For $0<s< p \leq q $, by H\"{o}lder inequality we get

\begin{equation}\label{3.11}
\begin{split}
 \int_{B_R} v^{-s }
    \leq & \Big( \int_{B_R}\big( v^{-s }\big)^{\frac{q }{s }}\Big)^{\frac{s }{q }}\Big( \int_{B_R}1^{\frac{q }{q- s}}\Big)^{\frac{q- s}{q }}\\
\lesssim & \Big(R^{n-q }\Big)^{\frac{s }{q }}R^{n\frac{q- s}{q }}\\
=& R^{n-s },
\end{split}
\end{equation}
where in the last second step we have used (\ref{3.3}) with $p \leq q<a+1$.
This implies (\ref{3.3}) also valid for $0<q< p $ and hence for all $0<q<a+1$.

To prove  (\ref{3.3}) for $q=a+1$, we need more careful calculations, since then the first term in (\ref{3.5}) vanishes and now we have

\begin{equation}\label{3.12}
\begin{split}
b\int \eta^\theta v^{-a-1} \lesssim \frac{1}{R}\int \eta^{\theta-1}v^{-a}|\nabla v|^{p-1},
\end{split}
\end{equation}
or that

\begin{equation}\label{3.13}
\begin{split}
\int_{B_R} v^{-a-1} \lesssim \frac{1}{R}\int_{B_{2R}} v^{-a}|\nabla v|^{p-1}.
\end{split}
\end{equation}
Next we estimate the right hand side of (\ref{3.13}). By H\"{o}lder inequality,

\begin{equation}\label{3.14}
\begin{split}
\,& \int_{B_{2R}} v^{-a}|\nabla v|^{p-1}\\
    =&\int_{B_{2R}} v^{(\varepsilon-a-1)\frac{p-1}{p}}|\nabla v|^{p-1}\cdot v^{-a-(\varepsilon-a-1)\frac{p-1}{p}}\\
\leq & \Big( \int_{B_{2R}} \big( v^{(\varepsilon-a-1)\frac{p-1}{p}}|\nabla v|^{p-1}\big)^{\frac{p}{p-1}}\Big)^{\frac{p-1}{p}}
       \Big( \int_{B_{2R}} \big( v^{-a-(\varepsilon-a-1)\frac{p-1}{p}}\big)^{p}\Big)^{\frac{1}{p}}\\
   =& \Big( \int_{B_{2R}}v^{\varepsilon-a-1}|\nabla v|^{p}\Big)^{\frac{p-1}{p}}
       \Big( \int_{B_{2R}}  v^{-a+(1-\varepsilon)(p-1)} \Big)^{\frac{1}{p}}
\end{split}
\end{equation}
Choosing $0<\varepsilon< \min\{1, \frac{n-p}{p}(p-1)\}$ implies $p<a+1-\varepsilon<a+1$ and $0<a-(1-\varepsilon)(p-1)<a+1$.
Then for the right hand side of (\ref{3.14}), we can use (\ref{3.3}) with $0<q<a+1$ and  (\ref{3.2})  to get

\begin{equation}\label{3.15}
\begin{split}
\,& \int_{B_{2R}} v^{-a}|\nabla v|^{p-1}\\
\lesssim & \Big( (2R)^{n+\varepsilon-a-1}\Big)^{\frac{p-1}{p}}
       \Big( (2R)^{n-a+(1-\varepsilon)(p-1)} \Big)^{\frac{1}{p}}\\
\lesssim & R^{n-a}.
\end{split}
\end{equation}
Submitting this into  (\ref{3.13}) we complete the proof of (\ref{3.3}) for $q=a+1$ as desired.\qed

\vspace{20pt}
$\mathbf{Proof\,\, of\,\, Theorem\,\, \ref{Thm1}}$ \qquad Now we give the proof of theorem \ref{Thm1}.

Let $\eta$ be the smooth cut-off functions as in  (\ref{3.1}) and take a constant $\theta>0$ big enough.
Replacing  $m$ with  $-m$ and $\varphi$ with $\eta^\theta$ in  (\ref{2.15}) we have

\begin{equation}\label{3.16}
\begin{split}
\,&\int_{\mathbb{R}^n}\eta^\theta v^{1-n}g^{-m} \mathbf{Tr}\{E^2\} -nm\int_{\mathbb{R}^n}\varphi v^{1-n}g^{-m-1}X^{j}E_{ij}E_i\\
\,&\qquad \qquad\qquad \qquad \qquad         \leq -\theta\int_{\mathbb{R}^n} \eta^{\theta-1} v^{1-n}g^{-m}X^jE_{ij}\eta_i .
\end{split}
\end{equation}

In the following, we will choose $m=\frac{p-1}{p}-\varepsilon_0>0$, with $\varepsilon_0>0$ small enough and depending only on $n,p$,
to deduce  (\ref{3.32}). Roughly speaking, first by using lemma \ref{lem-8} we will show the left hand side of  (\ref{3.16})
is nonnegative and obtain (\ref{3.20}).
Then, using corollary \ref{cor-1} to the right hand side of  (\ref{3.20}) we get (\ref{3.23}).
Finally, we will use lemma  \ref{lem3.1} to estimate the right hand side of  (\ref{3.23}) to deduce  (\ref{3.32}).

Notice that $g=av^{-1}|\nabla v|^p +bv^{-1}$, then we rewrite  the first term in (\ref{3.16}) as

\begin{equation}\label{3.17}
\begin{split}
\int_{\mathbb{R}^n}\eta^\theta v^{1-n}g^{-m} \mathbf{Tr}\{E^2\}
        = &\int_{\mathbb{R}^n}\eta^\theta v^{1-n}g^{-m-1}(av^{-1}|\nabla v|^p +bv^{-1}) \mathbf{Tr}\{E^2\} \\
        = & b\int_{\mathbb{R}^n}\eta^\theta v^{-n}g^{-m-1}\mathbf{Tr}\{E^2\}\\
        \,&\, +a\int_{\mathbb{R}^n}\eta^\theta v^{-n}g^{-m-1}|\nabla v|^p \mathbf{Tr}\{E^2\}.
\end{split}
\end{equation}
Therefor the left hand side of (\ref{3.16}) can be rewritten as

\begin{equation}\label{3.18}
\begin{split}
\, &\int_{\mathbb{R}^n}\eta^\theta v^{1-n}g^{-m} \mathbf{Tr}\{E^2\}-nm\int_{\mathbb{R}^n}\varphi v^{1-n}g^{-m-1}X^{j}E_{ij}E_i\\
 = & b\int_{\mathbb{R}^n}\eta^\theta v^{-n}g^{-m-1}\mathbf{Tr}\{E^2\}+(a-nm)\int_{\mathbb{R}^n}\eta^\theta v^{-n}g^{-m-1}|\nabla v|^p \mathbf{Tr}\{E^2\}\\
 \,&\quad +nm\int_{\mathbb{R}^n}\eta^\theta v^{-n}g^{-m-1}|\nabla v|^p \mathbf{Tr}\{E^2\}-nm\int_{\mathbb{R}^n}\varphi v^{1-n}g^{-m-1}X^{j}E_{ij}E_i\\
 = & b\int_{\mathbb{R}^n}\eta^\theta v^{-n}g^{-m-1}\mathbf{Tr}\{E^2\}+(a-nm)\int_{\mathbb{R}^n}\eta^\theta v^{-n}g^{-m-1}|\nabla v|^p \mathbf{Tr}\{E^2\}\\
 \,&\quad +nm\int_{\mathbb{R}^n}\eta^\theta v^{-n}g^{-m-1}|\nabla v|^p \big(\mathbf{Tr}\{E^2\}-\frac{v_{j}}{|\nabla v|}E_{ij}\cdot\frac{v_{k}}{|\nabla v|}E_{ki}\big) .
\end{split}
\end{equation}

By lemma \ref{lem-8}, the bracket pair in above last integral is nonnegative.
Now if we take $m=\frac{p-1}{p}-\varepsilon_0$ with $0<\varepsilon_0<\frac{p-1}{p}$, then

\begin{equation}\label{3.19}
\begin{split}
   \,&\int_{\mathbb{R}^n}\eta^\theta v^{1-n}g^{-m} \mathbf{Tr}\{E^2\}-nm\int_{\mathbb{R}^n}\varphi v^{1-n}g^{-m-1}X^{j}E_{ij}E_i\\
\geq & b\int_{\mathbb{R}^n}\eta^\theta v^{-n}g^{-m-1}\mathbf{Tr}\{E^2\}+n\varepsilon_0\int_{\mathbb{R}^n}\eta^\theta v^{-n}g^{-m-1}|\nabla v|^p \mathbf{Tr}\{E^2\}\\
  =  & \frac{p}{p-1}\varepsilon_0\Big[\frac{p-1}{p\varepsilon_0}b\int_{\mathbb{R}^n}\eta^\theta v^{-n}g^{-m-1}\mathbf{Tr}\{E^2\}
                                      +a\int_{\mathbb{R}^n}\eta^\theta v^{-n}g^{-m-1}|\nabla v|^p \mathbf{Tr}\{E^2\}\Big]\\
\geq & \frac{p}{p-1}\varepsilon_0\int_{\mathbb{R}^n}\eta^\theta v^{1-n}g^{-m} \mathbf{Tr}\{E^2\}.
\end{split}
\end{equation}
Together this with  (\ref{3.16}) yields

\begin{equation}\label{3.20}
 \frac{p}{p-1}\varepsilon_0\int_{\mathbb{R}^n}\eta^\theta v^{1-n}g^{-m} \mathbf{Tr}\{E^2\}\leq\int_{\mathbb{R}^n} \eta^{\theta-1} v^{1-n}g^{-m}X^jE_{ij}\eta_i .
\end{equation}

On the other hand, if we take $B$ with the entries $B_{ij}=\epsilon \eta^{-1}\eta_iX^j$ for $\epsilon>0$ small,
then using corollary \ref{cor-1} we get

\begin{equation}\label{3.21}
\eta^{-1}\eta_iX^jE_{ij} \leq \frac{c(p)}{\epsilon}\eta^{-2}\eta_i\eta_jX^iX^j+ \epsilon \mathbf{Tr}\{E^2\}.
\end{equation}
Plugging  this into (\ref{3.20})  we obtain

\begin{equation}\label{3.22}
(\frac{p}{p-1}\varepsilon_0-\epsilon)\int_{\mathbb{R}^n}\eta^\theta v^{1-n}g^{-m} \mathbf{Tr}\{E^2\}
\leq \frac{c(p)}{\epsilon}\int_{\mathbb{R}^n} \eta^{\theta-2} v^{1-n}g^{-m}\eta_i\eta_jX^iX^j .
\end{equation}
Since $|\eta_i\eta_j|\lesssim \frac{1}{R^2}$ and $|X^iX^j|\leq |\nabla v|^{2p-2}$,
taking $\epsilon= \frac{p}{2(p-1)}\varepsilon_0$  we arrive at

\begin{equation}\label{3.23}
\int_{\mathbb{R}^n}\eta^\theta v^{1-n}g^{-m} \mathbf{Tr}\{E^2\}\lesssim \frac{1}{R^2}\int_{\mathbb{R}^n} \eta^{\theta-2} v^{1-n}g^{-m}|\nabla v|^{2p-2}.
\end{equation}
For the term on the right hand side of (\ref{3.23}), by

$$g^{-m}=\frac{1}{\big(av^{-1}|\nabla v|^p +bv^{-1}\big)^{\frac{p-1}{p}-\varepsilon_0}}\leq\frac{1}{\big(av^{-1}|\nabla v|^p\big)^{\frac{p-1}{p}-\varepsilon_0}} $$

\noindent we have
\begin{equation}\label{3.24}
\begin{split}
\int_{\mathbb{R}^n} \eta^{\theta-2} v^{1-n}g^{-m}|\nabla v|^{2p-2}
\lesssim &\int_{\mathbb{R}^n} \eta^{\theta-2} v^{1-n}\big(v^{-1}|\nabla v|^p\big)^{-(\frac{p-1}{p}-\varepsilon_0)}|\nabla v|^{2p-2}\\
        =&\int_{\mathbb{R}^n} \eta^{\theta-2} v^{2-n-\frac{1}{p}+\varepsilon_0}|\nabla v|^{p-1+p\varepsilon_0}\\
        \lesssim & \int_{B_{2R}} v^{2-n-\frac{1}{p}+\varepsilon_0}|\nabla v|^{p-1+p\varepsilon_0}\\
        =&  \int_{B_{2R}} v^{(\varepsilon_0-a-1)\frac{p-1+p\varepsilon_0}{p}}|\nabla v|^{p-1+p\varepsilon_0}\cdot v^{\widetilde{p}},
\end{split}
\end{equation}
with $\widetilde{p}=2-n-\frac{1}{p}+\varepsilon_0- (\varepsilon_0-a-1)\frac{p-1+p\varepsilon_0}{p}$.

Recall  $0<\varepsilon_0<\frac{p-1}{p}$. Now we assume in addition $0<\varepsilon_0<\min\{\frac{p-1}{p}, \frac{1}{p} \}$.
 Then, we can use the  H\"{o}lder inequality with exponent
pair $(\frac{p}{p-1+p\varepsilon_0},\frac{p}{1-p\varepsilon_0})$ to the last term in (\ref{3.24}) to deduce

\begin{equation}\label{3.25}
\begin{split}
\int_{\mathbb{R}^n} \eta^{\theta-2} v^{1-n}g^{-m}|\nabla v|^{2p-2}
\lesssim &  \int_{B_{2R}} v^{(\varepsilon_0-a-1)\frac{p-1+p\varepsilon_0}{p}}|\nabla v|^{p-1+p\varepsilon_0}\cdot v^{\widetilde{p}}\\
\leq &\Big(\int_{B_{2R}} v^{\varepsilon_0-a-1}|\nabla v|^{p}\Big)^{\frac{p-1+p\varepsilon_0}{p}}
 \cdot\Big(\int_{B_{2R}}v^{q}\Big)^{\frac{1-p\varepsilon_0}{p}},
\end{split}
\end{equation}
where $q= \widetilde{p} \cdot\frac{p}{1-p\varepsilon_0}$.

Next, we will use lemma \ref{lem3.1} to estimate the last two factors in (\ref{3.25}).
Then we assume furthermore $0<\varepsilon_0<\min\{\frac{p-1}{p}, \frac{1}{p}, \frac{p-1}{p}(n-p)\}$.
So, for the first one, by $0<\varepsilon_0< \frac{p-1}{p}(n-p)=a-(p-1)$ we have $-a-1<\varepsilon_0-a-1<-p$.
Then we can use (\ref{3.2}) to obtain

\begin{equation}\label{3.26}
\Big(\int_{B_{2R}} v^{\varepsilon_0-a-1}|\nabla v|^{p}\Big)^{\frac{p-1+p\varepsilon_0}{p}}
\lesssim R^{(n+\varepsilon_0-a-1)\frac{p-1+p\varepsilon_0}{p}}.
\end{equation}
For the second one, first we observe

 $$q= \widetilde{p}\cdot \frac{p}{1-p\varepsilon_0}= \frac{1}{p}\big[ 3p^2-2(n+1)p+n\big]+\big[\frac{3p-(n+1)}{\frac{1}{p}-\varepsilon_0}+1\big] \varepsilon_0.$$

\noindent So we consider the mater in the following two cases:

$$(i)\quad  \frac{n+1}{3}\leq p<\frac{n+1+\sqrt{(n+1)^2-3n}}{3} $$
\noindent and
$$(ii)\quad \frac{n+1+\sqrt{(n+1)^2-3n}}{3}\leq p<n.$$

$Case (i).$ In this case, we see  $-a-1\leq \frac{1}{p}\big[ 3p^2-2(n+1)p+n\big]<0$ and $\big[\frac{3p-(n+1)}{\frac{1}{p}-\varepsilon_0}+1\big]>0$.
So if we choose $\varepsilon_0>0$ small enough, say
$0<\varepsilon_0<\min\{\frac{p-1}{p}, \frac{1}{2p}, \frac{p-1}{p}(n-p), \frac{-[3p^2-2(n+1)p+n]}{p\big[\frac{3p-(n+1)}{\frac{1}{2p}}+1\big]}\}$,
then  $-a-1<q<0$. Therefor, we can use (\ref{3.3}) to obtain

 \begin{equation}\label{3.27}
\Big(\int_{B_{2R}}v^{q}\Big)^{\frac{1-p\varepsilon_0}{p}}
\lesssim R^{(n+q)\frac{1-p\varepsilon_0}{p}}=R^{n\cdot\frac{1-p\varepsilon_0}{p}+\widetilde{p}}.
\end{equation}

Combining  (\ref{3.26}), (\ref{3.27}) with (\ref{3.25}) we have

 \begin{equation}\label{3.28}
\begin{split}
\int_{\mathbb{R}^n} \eta^{\theta-2} v^{1-n}g^{-m}|\nabla v|^{2p-2}
\lesssim &  R^{(n+\varepsilon_0-a-1)\frac{p-1+p\varepsilon_0}{p}}\cdot R^{n\cdot\frac{1-p\varepsilon_0}{p}+\widetilde{p}}\\
= & R^{2-\frac{1}{p}+ \varepsilon_0}.
\end{split}
\end{equation}

$Case (ii).$ In this case we have $q>0$. Therefor, by (\ref{2.2}) we get

 \begin{equation}\label{3.29}
\Big(\int_{B_{2R}}v^{q}\Big)^{\frac{1-p\varepsilon_0}{p}}
\lesssim R^{(n+q\cdot\frac{p}{p-1})\frac{1-p\varepsilon_0}{p}}.
\end{equation}
 Submitting  (\ref{3.26}) and  (\ref{3.29}) into (\ref{3.25}) we have

 \begin{equation}\label{3.30}
\begin{split}
\int_{\mathbb{R}^n} \eta^{\theta-2} v^{1-n}g^{-m}|\nabla v|^{2p-2}
\lesssim &  R^{(n+\varepsilon_0-a-1)\frac{p-1+p\varepsilon_0}{p}}\cdot R^{(n+q\cdot\frac{p}{p-1})\frac{1-p\varepsilon_0}{p}}\\
= & R^{n+(\varepsilon_0-a-1)\frac{p-1+p\varepsilon_0}{p}+\widetilde{p} +\frac{\widetilde{p}}{p-1} }\\
= &  R^{2-\frac{1}{p}+\frac{ 1}{p^2(p-1)}[3p^2-2(n+1)p+n]+ \frac{1}{p-1}\big( a+1+\frac{1}{p}-\varepsilon_0\big)\varepsilon_0 }\\
= &  R^{2-\frac{2p-1}{p^2(p-1)}(n-p)+ \frac{1}{p-1}\big( a+1+\frac{1}{p}-\varepsilon_0\big)\varepsilon_0 }.
\end{split}
\end{equation}

Now for any $\frac{n+1}{3}\leq p<n$, together  (\ref{3.23}) with  (\ref{3.28}) or (\ref{3.30}) we have

\begin{equation}\label{3.31}
\int_{\mathbb{R}^n}\eta^\theta v^{1-n}g^{-m} \mathbf{Tr}\{E^2\}\lesssim R^{-s(\varepsilon_0)},
\end{equation}
with  $s(\varepsilon_0)=\frac{1}{p}- \varepsilon_0$ in case (i)
or $s(\varepsilon_0)=\frac{2p-1}{p^2(p-1)}(n-p)- \frac{1}{p-1}\big( a+1+\frac{1}{p}-\varepsilon_0\big)\varepsilon_0$ in case (ii).
Finally if  we chose $\varepsilon_0>0$ small such that

$0<\varepsilon_0<\min\big{\{}\frac{p-1}{p}, \frac{1}{2p}, \frac{p-1}{p}(n-p), \frac{-[3p^2-2(n+1)p+n]}{p\big[\frac{3p-(n+1)}{\frac{1}{2p}}+1\big]}, \frac{2p-1}{p^2(a+2)}(n-p)\big{\}},$

\noindent then, in both cases,

\begin{equation}\label{3.32}
\int_{\mathbb{R}^n}\eta^\theta v^{1-n}g^{-m} \mathbf{Tr}\{E^2\}\lesssim R^{-s}
\end{equation}
for some constant $s>0$ depending only on $n,p$. Letting $R\rightarrow \infty$ in (\ref{3.32}) we deduce

\begin{equation}\label{3.33}
\int_{\mathbb{R}^n}  v^{1-n}g^{-m} \mathbf{Tr}\{E^2\}\leq 0,
\end{equation}
which implies $E=0$ $a.e$ in $\mathbb{R}^n$, especially $E\equiv 0$ in $\Omega^{c}_{cr}$. Thus

$$v=C_1+C_2|x-x_0|^{\frac{p}{p-1}}$$
\noindent for some $C_1,C_2>0$ and $x_0\in \mathbb{R}^n$ and hence $u(x)=U_{\lambda,x_0}$ (see e.g. \cite{CFR2020} or \cite{CMR2022} for more details).
The proof of Theorem \ref{Thm1} is completed. \qed

\vspace{0.6cm}
{\it Acknowledgement.}
The author  would like to express his gratitude to Prof. Xi-Nan Ma for helpful discussion on this topic and useful advice on the
manuscript. Research of the author  was supported by National Natural Science Foundation of China (grants 11861016 and 12141105).

\end{document}